\documentclass[a4paper,11pt]{article}
\usepackage{amsmath,amsfonts,latexsym}
\usepackage{graphicx}
\usepackage{epsfig}
\usepackage{epstopdf}
\nonstopmode

\numberwithin{equation}{section}
\numberwithin{table}{section}
\numberwithin{figure}{section}

\newtheorem{Theorem}{Theorem}[section]
\newtheorem{Proposition}{Proposition}[section]
\newtheorem{Lemma}{Lemma}[section]

\newtheorem{Definition}{Definition}[section]
\newtheorem{Remark}{Remark}[section]

\newenvironment{Proofc}[1]{\smallskip\par\noindent\textsc{#1}\quad}%
  {\hfill$\Box$\bigskip\par}
\newenvironment{Proof}{\begin{Proofc}{Proof}}{\end{Proofc}}
  {\bigskip\par}
\renewcommand{\epsilon}{\varepsilon}
\renewcommand{\phi}{\varphi}
\renewcommand{\a}{{\alpha}}
\renewcommand{\b}{{\beta}}
\newcommand{\g}{{\gamma}}
\renewcommand{\d}{{\delta}}

\renewcommand{\l}{\lambda}

\newcommand{\s}{{\sigma}}

\newcommand{\G}{{\Gamma}}

\newcommand{\Og}{\Omega}

\newcommand{\pd}{{\partial}}

\newcommand{\cC}{{\cal C}}
\newcommand{\cF}{{\cal F}}
\newcommand{\cG}{{\cal G}}

\newcommand{\Z}{{\mathbb Z}}

\newcommand{\R}{{\mathbb R}}
\newcommand{\N}{{\mathbb N}}

\newcommand{\be}{\begin{equation}}
\newcommand{\ee}{\end{equation}}

\newcommand{\lims}{ {\limsup\limits_{h\to 0}}^*}
\newcommand{\limi}{ {\liminf\limits_{h\to 0}}\phantom{\hskip -1pt}_*}
\begin{document}
\title{On the  approximation of the principal eigenvalue for a class of nonlinear     elliptic operators}
 \author{ Isabeau Birindelli\thanks{Dip. di Matematica, ``Sapienza" Universit\`a  di Roma,   P.le Aldo  Moro 2, 00185 Roma, Italy {\tt isabeau@mat.uniroma1.it}.}
 \and
 Fabio Camilli\thanks{Dip. di Scienze di Base e Applicate per l'Ingegneria,  ``Sapienza" Universit{\`a}  di Roma, via Scarpa 16,
 00161 Roma, Italy, ({\tt e-mail:camilli@dmmm.uniroma1.it})}\and Italo Capuzzo Dolcetta\thanks{Dip. di Matematica, ``Sapienza" Universit\`a  di Roma,   P.le Aldo  Moro 2, 00185 Roma, Italy {\tt capuzzo@mat.uniroma1.it}.}}
\date{version: \today}
\maketitle
 \begin{abstract}
     We present a finite difference method to compute the principal eigenvalue and the corresponding eigenfunction for a large class of second order elliptic operators including notably  linear operators in nondivergence form and fully nonlinear operators. \\
       The principal eigenvalue is computed by solving a finite-dimensional nonlinear  min-max optimization problem.
     We prove the convergence of the method and we discuss its  implementation. Some examples where the exact solution is explicitly known show the effectiveness of the method.
 \end{abstract}
\begin{description}
\item [{\bf MSC 2000}:] 35J60,  35P30, 65M06.
\item [{\bf Keywords}:]
Principal eigenvalue, nonlinear elliptic operators, finite difference schemes, convergence.
\end{description}
 \section{Introduction}\label{intro}
Consider the elliptic self-adjoint operator
\begin{equation}\label{Lself}
Lu(x)=\pd_i\left(a_{ij}(x)\pd_{j}u(x)\right),
\end{equation}
where $a_{ij}=a_{ji}$ are smooth functions in $ \Omega$, a smooth bounded open subset of $\R^n$, satisfying $a_{ij}\xi_i \xi_j\ge \alpha|\xi|^2$ for some $\alpha>0$.
It is well-known that the minimum value $\lambda_1$ in the Rayleigh-Ritz variational formula
\[
\lambda_1= \inf_{\phi\in H^1_0(\Omega), \phi\not\equiv 0} \frac{-  \int_{\Omega} \phi(x)\, L\phi(x) \,\,dx\,
}{ \|\phi\|^2_{L^2(\Omega)}}=\inf_{\phi\in H^1_0(\Omega), \phi\not\equiv 0} \frac{\int_{\Omega} a_{ij}(x)\pd_{j} \phi(x)\pd_i\phi(x) \,\,dx\,
}{ \|\phi\|^2_{L^2(\Omega)}}\]
 is attained at some function $w_1$ satisfying
\[
\left\{ \begin{array}{ll}
  Lw_1(x)+\lambda_1 w_1(x)=0   \qquad& x\in \Og,  \\
   w_1(x)=0 & x\in \partial\Og.
    \end{array}
\right.
\]
The number $\lambda_1$ is usually referred to as the principal eigenvalue of $L$ in $\Omega$ and $w_1$ is the corresponding principal eigenfunction.
For operators of the form \eqref{Lself} and also more general linear operator in divergence form there is a vast literature on computational methods for the principal eigenvalue, see for example \cite{BO}, \cite {B},  \cite {H}, \cite{W}.\par
General non-divergence type elliptic operators, namely
\begin{equation}\label{Lnonself}
Lu(x)=a_{ij}(x)\pd_{ij}u(x)+b_i(x)\pd_i u(x)+c(x)u
 \end{equation}
are not self-adjoint and the spectral theory is then much more involved: in particular, the Rayleigh-Ritz variational formula is not available anymore.
In the seminal paper \cite{DV2} by M.D. Donsker and S.R.S. Varadhan, a min-max formula for the principal eigenvalue of a class of elliptic operators $L$ including (\ref{Lnonself}) was proved,  namely
\begin{equation}\label{PE2intro}
     \l_{1}=-\inf_{\phi \in C^2(\Og), \phi>0}\;\sup_{x\in\Og}\frac{L\phi(x)}{\phi(x)}.
\end{equation}
In that papers other  representation formulas for $\l_{1}$ were also proposed in terms of large deviations   and of the average long run time behavior of the positive semigroup generated by $L$.
A  further crucial step in that direction is the paper \cite{BNV} by H. Berestycki, L. Nirenberg and S.R.S. Varadhan, where the validity of formula (\ref{PE2intro}) is proved under mild smoothness assumptions ($\Omega$ a bounded open set and $a_{ij}\in C^0(\Omega)$, $b_i$, $c\in L^\infty(\Omega)$). Moreover  it is proved that \eqref{PE2intro} is  equivalent to
\[
 \l_1:=\sup\{\l\in\R:\, \exists\, \phi>0 \;\text{ such that}\;L\phi+\l \phi\le 0\quad\text{in $\Og$}\}.\,
\]
Following this path of ideas, notions of principal eigenvalue for fully nonlinear uniformly elliptic operators of the form
$$F[u]= F(x, u(x), Du(x), D^2 u(x))$$
have been introduced and analyzed in \cite{A},  \cite{BCDPR}, \cite{BD}, \cite{BEQ},  \cite {IY}, \cite {L}.  A by now established definition of  principal eigenvalue is given by
\begin{equation}\label{PEC}
\l_1:=\sup\{\l\in\R:\, \exists\, \phi>0 \;\text{ such that}\;F[\phi]+\l \phi\le 0\quad\text{in $\Og$}\}\,
\end{equation}
where the inequality in \eqref{PEC} is intended in viscosity sense. It is possible to prove
under appropriate assumptions, see \eqref{Hyp_gen}-\eqref{Hyp_om}, that there exists a viscosity  solution $w_1$ of
\begin{equation}\label{PE}
\left\{ \begin{array}{ll}
  F[w_1]+ \l_1 w_1(x)=0   \qquad& x\in \Og,  \\
  w_1(x)=0 & x\in \partial\Og.
   \end{array}
\right.
\end{equation}
Moreover the characterization \eqref{PE2intro} still holds in this nonlinear setting.\par
As it is well-known, the principal eigenvalue plays a key role in several respects, both in the existence theory and in the qualitative analysis of  elliptic partial differential equations as well in applications to large deviations  \cite {A}, \cite {DV2}, bifurcation issues \cite {L}, ergodic and long run average cost problems in stochastic control \cite{BEN}. For linear non self-adjoint   operators and, a fortiori, for nonlinear ones the principal eigenvalue can be explicitly computed only in very special  cases,
see e.g. \cite{BL,Pu}, hence the importance to devise numerical algorithms for the problem. But, apart  some specific case (see \cite{BEM} for
the $p$-Laplace operator), approximation schemes and computational methods  are not available in the literature, at least at our present knowledge.\par
The aim of this paper is to define a numerical scheme for  the principal eigenvalue of nonlinear uniformly elliptic operators via a finite difference approximation of formula \eqref{PE2intro}.
More precisely, denoting by $\Z^n_h=h\Z^n$   the orthogonal lattice in $\R^n$ where $h>0$ is a discretization parameter, we consider a discrete operator $F_h$ acting on functions defined on a discrete subset $\Og_h\subset\Z^n_h$ of $\Omega$ and the corresponding approximated version of  \eqref{PE2intro}, namely
\begin{equation}\label{PE3intro}
     \l_{1,h}=-\inf_{\phi>0}\sup_{x\in\Og_h}\frac{F_h[\phi](x)}{\phi(x)}.
\end{equation}
As for the approximating operators $F_h$, we  consider a  specific class of finite difference schemes introduced in \cite{KT0}, \cite{KT1} since they  satisfy some useful properties for the convergence analysis.\par
We prove that if $F$ is   uniformly elliptic and satisfies in addition some quite natural further conditions, then it is possible to define a finite difference scheme $F_h$  such that   the discrete principal eigenvalues $\l_{1,h}$   and the associated discrete eigenfunctions $w_{1,h }$ converge  uniformly in $\Og$, as the mesh step $h$ is sent to $0$, respectively  to the principal eigenvalue $\l_1$  and to the corresponding eigenfunction $w_1$ for the original problem (\ref{PE}).
It is worth    pointing out that the proof of our main convergence result, Theorem \ref{main}, cannot rely on standard stability results    for fully nonlinear partial differential equations, see \cite{BS}, since the limit problem does not satisfy a   comparison principle (see Remark \ref{convergence} for details).\par
We  mention that our approach   is partially inspired by the paper \cite{GO} where a similar approximation scheme is proposed for the computation of effective Hamiltonians occurring in the homogenization of Hamilton-Jacobi equations which can be characterized by a formula somewhat similar to \eqref{PE2intro}.\par
In Section \ref{sect2} we introduce the main assumptions and we investigate some issues related to the Maximum Principle for discrete operators. In Section \ref{sect3} we study the approximation method for a class  of finite difference schemes and we prove the convergence of the scheme.  In Section \ref{sect4} we show that  under some additional structural assumptions on $F_h$ the inf-sup problem \eqref{PE3intro} can be transformed into a convex optimization problem on the nodes of the grid and we discuss its implementation. A few tests which show the efficiency of our method on some simple examples are reported in Section \ref{sect4} as well.
\section{The Maximum Principle for discrete operators}\label{sect2}
We start by fixing some notations and the  assumptions   on  the operator $F$.
Set $\Gamma=\Og\times\R\times\R^n\times S^n$ , where
$S^n$ denotes the linear space of real, symmetric $n\times n$ matrices. The function $F(x,z,p,r)$ is
assumed to be continuous on $\Gamma$ and locally
uniformly Lipschitz  continuous with respect to $z,p,r$ for each fixed $x\in \Og$. We will  also suppose
that the partial
derivatives $F_r$, $F_p$, $F_z$ satisfy the following
 structure conditions:
\begin{equation}\label{Hyp_gen}
 0<a I \le  F_r\le A I, \quad
 |F_p|\le \mu_1, \quad
 -\mu_0\le F_z\le 0.
\end{equation}
for some constants $a$, $A$, $\mu_0$, $\mu_1$.
A further condition is the positive homogeneity of degree $1$, that is
\begin{equation}\label{Hyp_om}
    F(x,tz,tp,tr)=tF(x,z,p,r)\qquad \forall t\geq 0.
\end{equation}
The principal eigenvalue of problem \eqref{PE} is defined by
\[
 \l_1=\sup\{\l:\, \exists\, \phi>0 \; { \rm such \,that} \;F[\phi]+\l\phi\le 0\ \mbox{in}\ \Og\},
\]
where the differential inequality $F[\phi]+\l\phi\le 0$ is  meant in the viscosity sense.
Under assumptions \eqref{Hyp_gen}-\eqref{Hyp_om}, there exists a viscosity solution of \eqref{PE} and the characterization \eqref{PE2intro}
of $\l_1$ holds (see \cite{BD}, \cite{BEQ}).

\begin{Remark}\label{negativeigenv}
It is possible to define
\[
 \l_1^-=\sup\{\l:\, \exists\, \phi<0 \; { \rm such \,that} \;F[\phi]+\l\phi\ge 0\ \mbox{in}\ \Og\}.
\]
When $F$ is not odd in its dependence on the Hessian, then in general $\l_1\neq\l_1^-$. Of course it is possible to see  $\l_1^-$ as $\l_1$ of some other operator. Hence we will only consider in this paper $\l_1$.
For example, for the extremal Pucci operators
${\mathcal  M}^+_{a,A}(D^2u):=\sup_{aI\leq B\leq AI}tr(AD^2u)$ and
${\mathcal  M}^-_{a,A}(D^2u):=\inf_{aI\leq B\leq AI}tr(AD^2u)$,
since ${\mathcal  M}^+_{a,A}(-M)=-{\mathcal  M}^-_{a,A}(M)$, the following
holds
$$\l_1^-({\mathcal  M}^+_{a,A})=\l_1({\mathcal  M}^-_{a,A}).$$

\end{Remark}
\begin{Remark}\label{cpositive}
The assumption $ F_z\le 0$, i.e. the monotonicity  of the differential operator in the zero-order term, could be removed.
Indeed  $\bar F:=F-c_0 z$,  with $c_0$ large, satisfies this assumption, moreover $\bar F$ and $F$ have the same principal
eigenfunction and the eigenvalues differ by $c_0$.
\end{Remark}
We now describe the discrete setting that we shall consider.
Given $h>0$, let $\Z^n_h=h\Z^n$ denote the orthogonal lattice in $\R^n$.
Let $F_h$  be a discrete operator acting on functions defined in
$\Og_h\subset\Z^n_h$. We shall consider an approximation of \eqref{PE} (which can be seen also as an eigenvalue problem for the discrete operator $F_h$). We look for  a number  $\l$ and a positive function  $w$ such that
\begin{equation}\label{PEd}
\left\{
\begin{array}{ll}
F_h(x,w (x),[w]_x)+\l w(x) =0 \qquad &x\in \Og_h,\\
 w(x)=0 & x\in \partial\Og_h,
\end{array}
\right.
\end{equation}
where
\begin{itemize}
\item $h>0$ is the discretization parameter ($h$ is meant to  tend to $0$),
\item $x\in \Og_h$ is the point where \eqref{PE}  is approximated,
\item $w$ is  a real valued mesh function  in $\Z^n_h$   meant to approximate the viscosity solution of \eqref{PE},
\item $[\cdot]_x$ represents the stencil of the scheme, i.e. the points in $\Og_h \backslash \{x\} $  where the value of $u$ is computed  for writing the scheme at the point $x$ (we  assume that $[w]_x$ is independent of      $w(y)$ for $|x-y|>Mh$ for some fixed $M\in\N$).
\end{itemize}
We denote by $\cC_h$ the space of the mesh functions defined on $\overline \Og_h$ and
 we introduce some basic assumptions for the scheme  $F_h$  (see \cite{KT0}, \cite{KT1}).
\begin{itemize}
\item[(i)] The operator $F_h$ is of positive type, i.e.
for all $x\in\Og_h$, $z,\tau\in\R$,  $u,\eta\in \cC_h$ satisfying $0\le \eta(y)\le \tau$ for each $y\in\Og_h$, then
\[F_h(x,z,[u+\eta]_x)\ge F_h(x,z,[u]_x)\ge F_h(x,z+\tau, [u+\eta]_x)\]
\item[(ii)] The operator $F_h$ is positively  homogeneous, i.e.
  for all $x\in\Og_h$, $z\in\R$,    $u\in \cC_h$ and $t\geq 0$, then
\[F_h(x,tz,[tu]_x)= t F_h(x,z,  [u]_x).\]
\item[(iii)]  The family of operators $\{F_h, 0<h\le h_0\}$, where $h_0$ is a positive constant,
is consistent  with the operator $F$ on the domain $\Og\subset\R^n$, i.e. for each $u\in C^2(\Og)$
 \[
\sup_{\Og_h} \left|F(x,u(x), Du(x),D^2u(x))-F_h(x,u(x),[u]_x)\right|\to 0\quad\text{as  $h\to 0$,}
 \]
 uniformly on compact subset  of $\Og$.
\end{itemize}

We study below some properties related to the maximum principle  and a comparison result for the operator $F_h$.
Let us start by the following definitions:
\begin{Definition}
A function $u\in \cC_h$ is a subsolution (respectively $v\in \cC_h$  is a supersolution) of
\begin{equation}\label{PEd2}
   F_h(x,u (x),[u]_x) =f(x) \qquad x\in\Og_h
\end{equation}
 if
 \begin{eqnarray*}
 &F_h(x,u (x),[u]_x) \ge f(x),\qquad x\in\Og_h\\
& \left(\text{respectively,  }   F_h(x,v (x),[v]_x) \le f(x),\qquad x\in\Og_h\right).
 \end{eqnarray*}
 \end{Definition}

\begin{Definition}\label{def_MP}
The Maximum Principle holds for the operator $F_h$ in $\Og_h$ if
\begin{equation}\label{MP}
\left\{\begin{array}{ll}
    F_h(x,u(x),[u]_x)\ge 0\qquad &\text{in $\Og_h$,}\\
    u\le 0&\text{on $\partial  \Og_h$,}
    \end{array}
    \right.
\end{equation}
implies $u\le 0$ in $\Og_h$.
\end{Definition}
\begin{Proposition}\label{prop_MP}
Assume that $F_h$ is of positive type and positive homogeneous and satisfies  either
\begin{equation}\label{C5_1}
  \begin{split}
 \text{for all $z\in\R$,  $u,\eta\in \cC_h$ satisfying $0\le  \eta(y)$  and $max_{y\in [\cdot]_x}\eta(y)>0$,}\\
\text{then }\quad F_h(x,z,[u+\eta]_x)>F_h(x,z,[u]_x)\qquad
\end{split}
\end{equation}
or
\begin{equation}\label{C5_2}
\begin{split}
   \text{for all $z,\tau\in\R$,  $u,\eta\in \cC_h$ satisfying $0\le \eta(y)\le \tau$ for each $y$, then}\\
F_h(x,z,[u]_x)\ge F_h(x,z+\tau, [u+\eta]_x)+c_0\tau\qquad
\end{split}
\end{equation}
for some positive constants  $c_0$.
Then   the Maximum Principle holds for the operator $F_h$ in $\Og_h$ .
\end{Proposition}
\begin{Proof}
Assume by contradiction that  $u$ satisfies \eqref{MP} and $M:=\max_{\overline \Og_h} u>0$. Let  $\bar x\in\Og_h$ be such that
$u(\bar x)=M$. Since $u\le0$ on $\partial\Og_h$, it is not restrictive to assume that there exists $y\in \Og_h$ such that $u(y)<u(\bar x)=M$.
%
Hence
\begin{align*}
    0&\le F_h(\bar x,u(\bar x), [u]_{\bar x})\le F_h(\bar x,u(\bar x)-M, [u-M]_{\bar x})
       \\
       & < F_h(\bar x,0,[0]_{\bar x})= 0,
\end{align*}
a contradiction. A similar  proof can be done with the assumption \eqref{C5_2}.

\end{Proof}
\begin{Remark}
The assumptions \eqref{C5_1} and \eqref{C5_2} correspond  to the  uniform ellipticity and, respectively, to the strict monotonicity of  the operator $F$  with respect  to the zero-order term.
\end{Remark}
The following proposition shows that, as it is known in the continuous case (see for example \cite{BNV,BD}), the validity of the Maximum Principle for subsolutions of  the operator $F_h$ is equivalent to the positivity  of the principal eigenvalue for $F_h$.
\begin{Proposition}\label{prop_MP2}
Assume that the scheme $F_h$  is of positive type and that it is positively
homogeneous.
Suppose that for $\l\in\R$, there exists a nonnegative grid function $\phi$ with $\phi>0$ in $\Og_h$  such that $F_h[\phi]+\l\phi\le 0$.
If, for $\tau<\l$,  the function $u$ satisfies
\[
\left\{\begin{array}{ll}
    F_h(x,u(x),[u]_x)+\tau u\ge 0\qquad &\text{in $\Og_h$}\\
    u\le 0&\text{on $\partial  \Og_h$,}
    \end{array}
    \right.
\]
then $u\le 0$ in $\Omega_h$, i.e. $F_{h,\tau}[\cdot]=F_h[\cdot]+\tau \cdot$ satisfies the Maximum Principle.
\end{Proposition}
\begin{Proof}
Suppose by contradiction that $\max_{\overline \Og_h}\{u\}>0$. Let $\phi$  as in the statement
 and set $L(\g)=\max_{ \Og_h}\{u-\g\phi\}$ (note that the maximum is taken only with respect to
 the  internal points). Then $L: [0,\infty)\to \R$ is continuous, decreasing, $L(0)>0$ and $L(\g)\to -\infty$ for $\g\to+\infty$.
 Hence there exists $\g'>0$  such that $L(\g')=0$. Moreover, since $u-\g'\phi\le 0$ on $\partial\Og$,
we also  have  $\max_{\overline \Og_h}\{u-\g'\phi\}=0$. Let $0<\g< \g'$ be such that
\begin{equation}\label{mp2}
\frac{\g}{\g'}\l>\tau
\end{equation}
and set $\psi=\g \phi$. Then $F_h[\psi]+\l \psi\le0$ and
$M=\max_{\overline \Og_h}\{u-\psi\}=(u-\psi)(\bar x)>0$ for some $\bar x\in \Og_h$. Hence $\psi(\bar x)+M=u(\bar x)$ and $\psi(x)+M\ge u(x)$.
Since $F_h$ is of positive type, it follows that
\begin{align*}
    F_h(x,\psi(\bar x),[\psi]_{\bar x})&\ge F_h((x,\psi(\bar x)+M,[\psi+M]_{\bar x})=F_h(x,u(\bar x),[\psi+M]_{\bar x})\\
                                                 &\ge F_h(x,u(\bar x),[u]_{\bar x}).
\end{align*}
Then
\begin{align*}
    \tau u(\bar x)\ge -F_h[u](\bar x)\ge -F_h[\psi](\bar x)\ge \l \psi(\bar x)
    =\l\g \phi(\bar x)\ge\l\frac{\g}{\g'}u(\bar x)
\end{align*}
and therefore a contradiction to \eqref{mp2}.
\end{Proof}
The following result gives    a comparison principle  for \eqref{PEd2}.
\begin{Proposition}\label{prop_comp}
Assume that $F_h$ is of positive type and it  satisfies either \eqref{C5_1} or \eqref{C5_2}. Let $u$ and $v$ be a subsolution
and respectively a supersolution  of \eqref{PEd2} such that $u\le v$ on $\pd \Og_h$. Then $u\le v$ in $\overline\Og_h$.
\end{Proposition}
\begin{Proof}
Suppose by contradiction that   $M:=\max_{\overline \Og_h} \{u-v\}>0$ and let  $\bar x\in\Og_h$ be such that
$u(\bar x)-v(\bar x)=M$. Hence  $v+M\ge u$ in $\Og_h$  and it is not restrictive to assume that $\max_{y_\in [\cdot]_{\bar x}}(v+M-u)>0$.
 It follows  that
\begin{align*}
f(\bar x)\le F_h(\bar x, u(\bar x), [u]_{\bar x})&=F_h(\bar x, v(\bar x)+M, [u]_{\bar x})
< F_h(\bar x, v(\bar x)+M, [v+M]_{\bar x})\\
&\le F_h(\bar x, v(\bar x), [v]_{\bar x})\le f(\bar x)
\end{align*}
and therefore a contradiction.
  A similar  proof can be carried on under assumption \eqref{C5_2}.
\end{Proof}


\section{Approximation  of the principal eigenvalue}
\label{sect3}
In this section we  consider a specific class of finite difference schemes introduced in
\cite{KT1}. These schemes satisfy certain
pointwise estimates which are  the discrete analogues of those valid for a general
class of fully nonlinear, uniformly elliptic equations.\\
We  assume that for all $x\in\Z^n_h$, the stencil $[\cdot]_x$ of the scheme  is given
by  $x+h Y$ where $Y=\{y_1,\dots,y_k\}\subset \Z^n$ is a finite set containing all the
vectors of the canonical basis of $\R^n$.
Then  we consider a discrete operator $F_h$  in \eqref{PEd}  given by a finite difference scheme  written in the form
\begin{equation}\label{finitediff}
F_h[u]={\cF} (x, u, \d_h u,\d^2_h u),
\end{equation}
where $\cF:\R^n\times\R\times\R^Y\times\R^Y\to \R$ and  for $y\in Y$, $u\in \cC_h$
\begin{align*}
  \delta_{h,y}^\pm u(x)&=\pm\frac{u(x\pm h y)-u(x)}{h|y|},\\
  \delta_{h,y} u(x)&=\frac{1}{2}\{ \delta_{h,y}^+ u(x)+ \delta_{h,y}^- u(x)\}=\frac{u(x+hy)-u(x-hy)}{2h|y|},\\
  \delta^2_{h,y} u(x)&=\delta_{h,y}^+\delta_{h,y}^- u(x)=\frac{u(x+hy)+u(x-hy)-2u(x) }{h^2|y|^2},\\
  \d_h u&=\{\d_{h,y} u:\, y\in Y\}, \quad
  \d^2_{h} u=\{\d^2_{h,y} u:\, y\in Y\}.
\end{align*}
Set  $\tilde\Gamma:=\R^n\times\R\times\R^Y\times\R^Y$ and denote by $(x,z,q,s)$ the generic points in $\tilde \Gamma$.
The operator $F_h$ given by \eqref{finitediff} is of positive type if
\begin{align}
  & \frac{\partial\cF}{\partial s_y}-\frac{|hy|}{2}\left|\frac{\partial\cF }{\partial q_y}\right|\ge 0\quad \forall y\in Y,\label{P1_scheme}\\[4pt]
  &\frac{\partial\cF}{\partial z}\le 0, \label{P2_scheme}
\end{align}
 and positively  homogeneous if
\[
  \cF (x,tz,tq,ts)=t \cF (x,z,q,s) \qquad \forall t\geq 0.
\]
 Moreover  if $F$ in \eqref{PE} satisfies the assumptions \eqref{Hyp_gen}, then it is always possible to find a scheme of type \eqref{finitediff} which is consistent with
$F$ and which, besides \eqref{P1_scheme}-\eqref{P2_scheme},   satisfies for all $y\in Y$,  the bounds
\begin{align}
    \frac{\partial\cF}{\partial s_y}-\frac{|hy|}{2}\left|\frac{\partial\cF }{\partial q_y}\right|\ge \a_0,\qquad
   \frac{\partial\cF }{\partial s_y}\le a_0, \qquad \left|\frac{\partial\cF}{\partial q_y}\right|\le b_0 \label{P3_scheme}
\end{align}
where  $\a_0$, $a_0$, $b_0$ are constants depending on $a$, $A$, $\mu_0$, $\mu_1$ in \eqref{Hyp_gen} (see     \cite{KT0}, \cite{KT1}). Note that  in particular \eqref{P3_scheme}  implies \eqref{C5_1}.\\
We recall some important properties of the previous scheme (for the proof   we refer to \cite{KT1})
\begin{Proposition}\label{prop_wellposed}
Assume \eqref{P1_scheme}-\eqref{P3_scheme} and let   $f$, $g$ be  two given mesh functions. Then for every $h>0$ sufficiently small  there exists a unique   solution $u_h:\Og_h\to\R$  to the Dirichlet problem
\begin{equation}\label{Dird}
   \left\{ \begin{array}{ll}
     F_h(x,u(x),[u]_x)=f\quad& x\in \Og_h, \\
     u=g &  x\in \pd\Og_h.
    \end{array}
    \right.
\end{equation}
\end{Proposition}
\begin{Proposition}\label{Prop_ABP}
Assume \eqref{P1_scheme}-\eqref{P3_scheme} and let  $u_h$ be a  subsolution of \eqref{Dird}. Then
\begin{equation}\label{ABP}
    \max_{\Og_h} u_h\le \max_{\pd \Og_h} g+ \frac{C}{\a_0}\left\{\sum_{x\in\Og_h}h^n|f(x)|^{n}\right\}^{\frac{1}{n}},
\end{equation}
where the constant $C$ is independent of $h$. Moreover if $u_h$ is a solution of \eqref{Dird}, then for any $x,y\in\Og_h$
\begin{equation}\label{Holder}
   |u_h(x)-u_h(y)|\le C \frac{|x-y|^\d}{R}\left(\max_{B_R^h} u_h+\frac{R}{\a_0}\left\{\sum_{x\in\Og_h}h^n|f(x)|^{n}\right\}^{\frac{1}{n}}\right),
\end{equation}
where $R=\min\{\mathrm{dist}(x,\pd \Og_h),\mathrm{dist}(x,\pd \Og_h)\}$, $B_R^h=B(0,R)\cap\Og_h$,  $\d$ and  $C$ are positive constants  independent of $h$.
\end{Proposition}
We give an example of a scheme of the form \eqref{finitediff}. Consider the Hamilton-Jacobi-Bellman operator \[
F(x,u,D u(x), D^2u(x))=\sup_{\a\in A}\inf_{\b\in B} L^{\a\b}u(x)
\]
where
\begin{equation}\label{op_L}
L^{\a\b}u(x)= a_{ij}^{\a\b}(x)D_{ij}u+b_i^{\a\b}(x)D_iu(x)+c^{\a\b}(x)u(x).
\end{equation}
It is always possible to rewrite  the operator $L^{\a\b}$ in \eqref{op_L}
  in the following form (see \cite{KT1})
\[
\overline L^{\a\b}  u(x)= \bar  a_{k}^{\a\b}(x)D^2_{y_k}u+\bar b_k^{\a\b}(x)D_{y_k}u(x)+\bar c^{\a\b}(x)u(x)
\]
where $D_{y_k}u =\langle Du,y_k\rangle$ and  $Y=\{y_1,\dots,y_k\}\subset \Z^n$ is a finite set containing all the vectors of the canonical basis
in $\R^n$. Moreover the coefficients $\bar  a_{k}^{\a\b}$, $\bar b_k^{\a\b}$ and $\bar c^{\a\b}$ satisfy the same properties  of $ a_{ij}^{\a\b}$, $ b_{ij}^{\a\b}$ and $c^{\a\b}$.
Then we consider
\begin{equation}\label{HJBd}
F_h[u](x):=\sup_{\a\in A}\inf_{\b\in B} L_h^{\a\b}u(x)
\end{equation}
where
\begin{equation}\label{op_Ld}
L_h^{\a\b}u(x)=  \bar  a_{k}^{\a\b}(x) \delta^2_{h,y_k} u(x)+\bar b_k^{\a\b}(x) \delta_{h,y_k}u(x)+\bar c^{\a\b}(x)u(x).
\end{equation}
For $x\in\R$ with $Y=\{1\}$ the previous scheme  reads as
\[
\begin{split}
   \sup_{\a\in A}\inf_{\b\in B}\Big\{&  a^{\a\b}(x)\frac{u(x+h)+u(x-h)-2u(x)}{h^2}+  b^{\a\b}(x)\frac{u(x+h)-u(x-h)}{2h}+\\
   &+  c^{\a\b}(x)u(x)\Big\}=0.
   \end{split}
\]


\subsection{The linear case}
In this part we assume that the operator $F$ in \eqref{PE} is linear, i.e. $F[u]=Lu$ with
\[
     Lu=a_{ij}(x)D_{ij}u+b_i(x)D_iu(x)+c(x)u(x)
\]
and we consider a scheme  defined  as in \eqref{HJBd}--\eqref{op_Ld}, obviously without the
dependence on $\a$, $\b$.

\begin{Proposition}\label{KR}
Under the assumption  \eqref{P3_scheme} the eigenvalue problem \eqref{PEd} has a simple eigenvalue $\l_{1,h}\in\R$ which corresponds to a positive eigenfunction. The other eigenvalues correspond to sign changing eigenfunctions.
\end{Proposition}
\begin{Proof}
Choose $\xi>0$ large enough so that $c(x)-\xi<0$ and set
$$L_{h,\xi}(x,t,[ u]_x)=L_h (x,t,[u]_x)-\xi t.$$
Let $K$ be the positive cone of the nonnegative grid functions in $\cC_h$. For a given grid  function $f$, by Proposition  \ref{prop_wellposed} and Proposition \ref{prop_comp} there exists a unique
solution $u\in \cC_h$  to
\[
   \left\{
   \begin{array}{ll}
    L_{h,\xi}(y,u(y),[u]_y)+f=0\qquad &\text{in $\Og_h$,}\\
    u= 0&\text{on $\partial  \Og_h$.}
   \end{array}
 \right.
\]
Since $\cC_h$ is a finite dimensional space it follows that $T:\cC_h\to \cC_h$ defined by $Tf=u$ is a compact linear operator.
Moreover, if $f\ge 0$,  then by Proposition \ref{prop_MP}  $u\ge 0$ and if $f\in K\setminus \{0\}$, $u=Tf>0$. \\Therefore,
by the Krein-Rutman theorem \cite{KR}, $r(T)$ the spectral radius of $T$ is a simple real  eigenvalue $r(T)>0$ with a positive  eigenfunction $u$
such that $Tu=r(T)u$. Hence for $\l_{1,h}=r(T)^{-1}-\xi$,  $w_1=Tu$ satisfies
\[
 \left\{\begin{array}{ll}
    L_h (x,w_1(x),[w_1]_x)  +\l_{1,h}w_1=0\qquad &\text{in $\Og_h$,}\\
    w_1= 0&\text{on $\partial  \Og_h$}.
   \end{array}
 \right.
\]
\end{Proof}
The following characterization of $\l_{1,h}$ is a simple consequence of Proposition \ref{prop_MP2}.
\begin{Proposition}
We have
\begin{equation}\label{PEd_char1}
    \l_{1,h}=\sup\left\{\l:\,\text {$\exists$  $\phi> 0$ s.t. $L_h[\phi]+\l \phi\le 0\ \mbox{in}\ \Og $}\right\},
\end{equation}
or, equivalently,
\begin{equation}\label{PEd_char2}
     \l_{1,h}=-\inf_{\phi>0}\sup_{x\in\Og_h}\left\{\frac{L_h[\phi](x)}{\phi(x)}\right\}.
\end{equation}
\end{Proposition}
\begin{Proof}
Denote by $\bar \l$ the right hand side of \eqref{PEd_char1}. Clearly $\l_{1,h}\le \bar\l$. If $\l_{1,h}< \bar\l$ then there exist $\mu\in (\l_{1,h}, \bar\l)$ and $\phi>0$ such that $L_h[\phi]+\mu \phi\le 0$.
A contradiction follows immediately by Proposition \ref{prop_MP2} since the eigenfunction corresponding to $\l_{1,h}$ is positive.
Hence  we have \eqref{PEd_char1}.\\
Let $\phi>0$ such that $L_h[\phi](x)+\l \phi(x)\le 0$  for $x\in\Og_h$. Hence
\[\l\le \inf_{\Og_h}\left\{-\frac{L_h[\phi]}{\phi}\right\}= -\sup_{\Og_h}\left\{\frac{L_h[\phi]}{\phi}\right\}.\]
Consequently
\[\l_{1,h}=\sup_{\phi>0}\left(-\sup_{x\in\Og_h}\left\{\frac{L_h[\phi](x)}{\phi(x)}\right\}\right)=-\inf_{\phi>0}\sup_{x\in\Og_h}\left\{\frac{L_h[\phi](x)}{\phi(x)}\right\}.\]
\end{Proof}
We give next an upper bound for $\l_{1,h}$  (compare with the corresponding estimate for $\l_1$ in \cite{BNV}, Lemma 1.1).
\begin{Lemma}
Let $n=1$ and assume that $B_R=\{|x|<R\}$ lies in $\Og$ with $R\le 1$. Then
\[
    \l_{1,h}(\Og_h)\le \frac{C}{R^2}
\]
\end{Lemma}
\begin{Proof}
Given  the linear operator
\[Lu=a(x)u''+b(x)u'(x)+c(x)u(x),\]
let $\g_0$, $\G_0$, $b$  be   positive constants  such that $\g_0\le a (x)\le \G_0$ and $|b(x)|,\, |c(x)|\le b$ in $\Og$.
Let $r=R/2$ and assume for simplicity that $r=Nh$ for some $N\in \N$. Set  $B_r=\{|x|<r\}$ and consider the grid  function
\[
    \s_i=(r^2-|ih|^2)^2 \qquad i=-N+1,\dots\,N-1
\]
Then for $i=-N+1,\dots, N-1$ we have
\begin{align*}
  &\frac{\s_{i+1}-\s_{i-1}}{2h}=-4hi(r^2-|ih|^2)+4h^3i\\
   &\frac{\s_{i+1}+\s_{i-1}-2\s_i}{h^2}=- 4(r^2-|ih|^2)+2h^2(2i^2+1).
\end{align*}
Denote by $a_i$, $b_i$ and $c_i$  the  coefficients of the linear operator computed at the point $x=ih$.
Since $h^2/(r^2-|ih|^2)\le 1$ it follows that
\begin{equation}\label{st_auto1}
\begin{split}
  -\frac{ L_h[\s](ih)}{4\s_i}&\le   \frac{ a_i}{(r^2-|ih|^2)}-\frac{a_i |hi|^2}{(r^2-|ih|^2)^2}
   +|b_i|\frac{2r}{(r^2-|ih|^2)}+\frac{c_i}{4}\\
   &    \le \frac{ \G_0+br}{(r^2-|ih|^2)}-  \frac{\g_0|hi|^2}{(r^2-|ih|^2)^2}+\frac{b}{4}.
\end{split}
\end{equation}
If
\[ |ih|^2(\g_0+\G_0+br)>r^2(\G_0+br),\]
then the second term in \eqref{st_auto1} dominates the first one and therefore
\begin{equation}\label{st_auto2}
-\frac{ L_h[\s](ih)}{4\s_i}\le \frac{b}{4}.
\end{equation}
In the remaining part of $B_r$,
\begin{equation}\label{st_auto3}
   -\frac{ L_h[\s](ih)}{4\s_i}\le \frac{ \G_0+br}{(r^2-|ih|^2)}+\frac{b}{4}\le \frac{b}{4}+\frac{1}{\g_0r^2}(\G_0+br)(\g_0+\G_0+br).
\end{equation}
By \eqref{st_auto2} and \eqref{st_auto3}, we get
\[
 \sup_{B_r}\left(-\frac{L_h[\s](ih)}{\s_i}\right)\le\frac{C}{R^2}\quad\text{for $i=-N+1,\dots, N-1$.}
\]
To conclude the proof, we show that  if for some positive function $\phi$ and $\l\in\R$, $L_h[\phi]+\l\phi\le 0$, then $\l\le \sup_{B_r}(-L_h[\s]/\s)$.
For this purpose, assume that $ \l> \sup_{B_r}(-L_h[\s]/\s):=\tau$; then $ L_h[\s]+\tau \s\ge 0$ in $B_r$ and  $\s=0$ on $\partial B_r$, while $L_h[\phi]+\tau \phi\le 0$. Hence by Proposition \ref{prop_MP2},
it follows $\s\le 0$ in $B_r$, a contradiction, and therefore $\l\le \tau$.
\end{Proof}
%
\subsection{The nonlinear case}
We consider now a general discrete operator $F_h$ given by \eqref{finitediff} and we study the
corresponding eigenvalue problem \eqref{PEd}. In analogy with formula \eqref{PEd_char1}, we define
\begin{equation}\label{PEd1}
    \l_{1,h}=\sup\left\{\l:\,\text {$\exists$  $\phi> 0$ \;{\rm such that} $F_h[\phi]+\l \phi\le 0$}\right\}
\end{equation}
We  prove  for each $h$ the existence  of a pair  $(\l_{1,h}, w_{1,h})$  satisfying   \eqref{PEd} with $w_{1,h}>0$ in $\Og_h$.
\begin{Proposition}\label{PEd:prop_exlower}
Assume that $F_h$ satisfies \eqref{P3_scheme}, $f\le 0$ and $\l< \l_{1,h}$.
Then there exists a nonnegative solution to
\begin{equation}\label{PEd3}
\left\{
\begin{array}{ll}
F_h(x,u (x),[u]_x)+ \l u(x) =f(x) \qquad &x\in \Og_h,\\
 u(x)=0 & x\in \partial\Og_h.
\end{array}
\right.
\end{equation}
\end{Proposition}
\begin{Proof}
We can assume $\l\ge 0$, since for  $\l<0$,   $F_h[u]+\l u$ satisfies  \eqref{C5_2} and therefore
by Propositions \ref{prop_wellposed} and \eqref{prop_comp} there exists a unique solution to problem \eqref{PEd3}.\\
Let us define by induction a sequence $u_n$ by setting $u_1\equiv 0$ and, for $n\ge1$ we consider the equation:
\begin{equation}\label{PEd4}
\left\{
\begin{array}{ll}
F_h(x,u_{n+1} (x),[u_{n+1}]_x) =f(x)-\l u_{n},\qquad &x\in \Og_h,\\
 u_{n+1}(x)=0 & x\in \partial\Og_h.
\end{array}
\right.
\end{equation}
\\
For any $n\in N$ there exists a non negative solution $u_{n+1}$ to \eqref{PEd4}. For $n=1$, existence follows by  Proposition \ref{prop_wellposed}. Moreover since $u_1\equiv 0$ is a subsolution to \eqref{PEd4}, by Proposition \ref{prop_comp} we get $u_2\ge 0$. The existence of a non negative solution at the $(n+1)$-step is proved in a similar way; moreover the solution is  non negative  since $f-\l u_{n}\le 0$. \\
We claim now that, for any $n\geq 1$, $u_n\leq u_{n+1}$. For $n=1$ the claim is trivially true since $u_2\ge 0$. Assume then by induction that $u_n\ge  u_{n-1} $.  Since $f(x)-\l u_{n}\le f(x)-\l u_{n-1}$ it follows that  $u_n$ is a subsolution of \eqref{PEd4}.
By  Proposition \ref{prop_comp}, we get that $u_n\le u_{n+1}$. \\
Let us show now that the sequence $u_n$ is bounded. Assume by contradiction that it is false and set $\overline u_n=u_n/|u_n|_\infty$. Then, by positive homogeneity, $\overline u_n$ is a solution of
\[F_h(x,\overline u_{n+1} (x),[\overline u_{n+1}]_x) =\frac{f(x)}{|u_{n+1}|_\infty}-\l \frac{u_{n}}{|u_{n+1}|_\infty},\qquad  x\in \Og_h.\]
Since the sequence $\overline u_n$ is bounded, then  up to a subsequence it converges to a function $\overline u $, while $u_{n}/|u_{n+1}|_\infty$ converges to $k \overline u$ where $k=\lim_{n\to \infty}  |u_{n}|_\infty/|u_{n+1}|_\infty\le 1$. Hence
$\overline u \ge 0$, $|\overline u |_\infty=1$, $\overline u =0$ on $\pd \Og_h$ and
\[F_h(x,\overline u (x),[\overline u ]_x) +k\l \overline u =0,\qquad  x\in \Og_h.\]
Since $0\le k\l\le \l$ and using the fact that for  $\l<\l_{1,h}$  there exists by definition $\phi>0$ such that $F_h[\phi]+\l\phi\ge 0$ in $\Og_h$, we get a contradiction to
Proposition \ref{prop_MP2}. Hence the sequence $u_n$ is bounded, and being in addition monotone, it converges pointwise to a function $u$ which solves \eqref{PEd3}.
\end{Proof}
The next result shows that $\l_{1,h}$ is indeed an eigenvalue for the approximated operator $F_h$.
\begin{Theorem}\label{PEd:prop_ex}
Assume that $F_h$ satisfies \eqref{P3_scheme}. Then there exists  $w_{1,h}>0$
in $\Og_h$ satisfying
\begin{equation}\label{PEd8}
\left\{
\begin{array}{ll}
F_h(x,w_{1,h}(x),[w_{1,h}]_x)+ \l_{1,h}w_{1,h}(x) =0 \qquad &x\in \Og_h,\\
w_{1,h}=0 & x\in \partial\Og_h.
\end{array}
\right.
\end{equation}
Moreover  that the characterization \eqref{PEd_char2} is still valid for the  nonlinear operator $F_h$.
\end{Theorem}
\begin{Proof}
Let $\l_n$ be an increasing sequence converging to $\l_{1,h}$. By Proposition \ref{PEd:prop_exlower} there exists a positive  solution $u_n$   of
\[
\left\{
\begin{array}{ll}
F_h(x,u_{n} (x),[u_{n}]_x)+\l_n u_{n} =-1,\qquad &x\in \Og_h,\\
 u_{n}(x)=0 & x\in \partial\Og_h.
\end{array}
\right.
\]
We claim that $u_n$ is not bounded. Assume by contradiction that $u_n$ is bounded so that, up to a subsequence, $u_n$ converges to
a function $u>0$ which solves
\[
\left\{
\begin{array}{ll}
F_h(x,u (x),[u]_x)+\l_{1,h} u =-1,\qquad &x\in \Og_h,\\
 u(x)=0 & x\in \partial\Og_h.
\end{array}
\right.
\]
Then, for $\epsilon>0$ small enough, $u$ satisfies
\[F_h(x,u (x),[u]_x)+(\l_{1,h}+\epsilon) u =-1+\epsilon u\le 0\]
which gives a contradiction to the definition \eqref{PEd1}. Hence $|u_n|_\infty\to \infty$. \\
Define now $w_n=u_n/|u_n|_\infty$ that solves
\[\left\{
\begin{array}{ll}
F_h(x,w_{n} (x),[w_{n}]_x)+\l_n w_{n} =-\frac{1}{|u_n|_\infty},\qquad &x\in \Og_h,\\
 w_{n}(x)=0 & x\in \partial\Og_h.
\end{array}
\right.\]
Then, up to a subsequence, $w_n$ converges to a bounded function $w_{1,h}$ which has norm 1 and which satisfies \eqref{PEd8}, so that $w_{1,h}>0$. \\
It is immediate that  \eqref{PEd_char2} is still valid for $F_h$.
\end{Proof}

\begin{Remark}\label{convergence}
There is a huge literature about the approximation of viscosity solutions of  first and second order PDEs. In this framework a well established technique to prove the   convergence of a numerical scheme is the Barles-Souganidis'method \cite{BS}: besides some natural properties of the scheme (stability, consistency, monotonicity), a key ingredient for  this technique is a \emph{strong comparison result} for the continuous problem which allows to show that a subsolution  is always lower than or equal to a supersolution. The comparison principle implies in particular that there is at most one viscosity solution of the problem. But it is immediate that  \eqref{PE} cannot   satisfy a comparison principle since $w\equiv 0$ and the principal eigenfunction $w_1$ are two distinct  solutions of the problem, hence the convergence proof cannot rely on the Barles-Souganidis'method and it needs a different argument.
\end{Remark}

We now discuss the convergence of the discrete principal eigenvalue
$\l_{1,h}$ to the continuous one defined by \eqref{PE2intro}.
We recall the definition of weak limits in viscosity sense (see \cite{BS})
\begin{align*}
 \lims u_h(x):=\lim_{h\to 0^+}\sup\{u_\d(y):\,|x-y|\le h,\, \d\le h\},\\
 \limi u_h(x):=\lim_{h\to 0^+}\inf\{u_\d(y):\,|x-y|\le h,\, \d\le h\}.
\end{align*}
\begin{Theorem}\label{main}
Assume \eqref{Hyp_gen}-\eqref{Hyp_om}, \eqref{P1_scheme}-\eqref{P3_scheme} and that $F_h$  is consistent with $F$.
Let $(\l_{1,h}, w_{1,h})$ be the sequence of the discrete  eigenvalues and of the corresponding eigenfunctions, solutions of \eqref{PEd}.
Then
$\l_{1,h}\to \l_{1}$ and $w_{1,h}\to w_{1}$ uniformly in $\overline \Og$ as $h\to 0$, where $\l_1$ and $w_1$ are  respectively the principal eigenvalue and a corresponding eigenfunction associated to $F$.
\end{Theorem}
\begin{Proof}\\
By the positive homogeneity of the scheme, it is not restrictive to assume that $\max_{\Og_h}\{w_{1,h}\}=1$, hence the sequence $w_{1,h}$ is bounded.
We first prove that
\begin{equation}\label{conv1}
     \liminf_{h\to 0}\l_{1,h}\ge \l_1.
\end{equation}
Assume by contradiction that
$\liminf_{h\to 0}\l_{1,h}=\tau$
for some some $\tau<\l_1$.
Consider  a subsequence, still denoted by $\l_{1,h}$, such that $\lim_{h\to 0} \l_{1,h}=\tau$.
Set $\overline w= \lims w_{1,h}$. By   standard stability results in viscosity solution theory, see \cite{BS}, $\overline w$ satisfies in viscosity sense
\begin{equation}\label{conv1a}
F[\overline w]+\tau \overline w\ge 0\qquad  \text{in $\Og_h$,}
\end{equation}
and
\begin{equation}\label{conv1b}
 \max_{\Og} \overline w = 1.
\end{equation}
 Let  $\eta>0 $ be such that for $h$ sufficiently small, $\l_{1,h}\le \tau+\eta$. Hence
\begin{align*}
    F_h[w_{1,h}]=-\l_{1,h}w_{1,h}\ge -\tau-\eta, \qquad x\in \Og_h.
\end{align*}
Set $f=-\tau-\eta$, $g\equiv 0$ and let $u_h$  be the corresponding
solution of \eqref{Dird}, while $u$ is the solution of
\[
   \left\{ \begin{array}{ll}
     F (x,u(x),Du, D^2u)=-\tau-\eta\quad& x\in \Og,  \\
     u=0 &  x\in \pd\Og.
    \end{array}
    \right.
\]
Then by Proposition \ref{prop_comp} and the consistency of the scheme   for $h$ sufficiently small
\begin{equation}\label{conv1d}
   0\le w_{1,h}\le u_h\le u+o(1)\qquad \text{in $\overline\Og_h$}
\end{equation}
and therefore
\begin{equation}\label{conv1c}
\overline w=0\quad  \text{on $\pd\Og$}.
\end{equation}
By \eqref{conv1a}, \eqref{conv1b} and \eqref{conv1c} we get a contradiction to the maximum principle for the operator $F$ (see \cite{BD}, \cite{BEQ}) and therefore \eqref{conv1}.\\
We now  prove that
\begin{equation}\label{conv2}
     \limsup_{h\to 0}\l_{1,h}\le \l_1.
\end{equation}
Assume by contradiction that there exists $\eta >0$ such that
 \[\bar \l:=\limsup_{h\to 0}\l_{1,h}\ge \l_1+\eta.\]
 We consider a subsequence,
still denoted by $\l_{1,h}$, such that $\lim_{h\to 0}\l_{1,h}=\bar \l$ and we set
 $\underline w= \limi w_{1,h}$. By  standard stability results
$\underline w$ satisfies $0\le \underline w\le 1$ and
\begin{equation}\label{conv3}
\left\{
   \begin{array}{ll}
    F[\underline w ]+(\l_1+\eta) \underline w\le 0\qquad &\mbox{ in } \Og,\\
    \underline w= 0&\mbox{ on }\ \partial  \Og,
    \end{array}
\right.
\end{equation}
in viscosity sense. Let    $x_h\in\Og_h$ be a   sequence such that $x_h\to x_0\in\overline\Og$ and  $w_{1,h}(x_h)=1$ for all $h>0$. By \eqref{conv1d}, $x_0\in\Og$.  We claim that
\begin{equation}\label{conv6}
\underline w(x_0)>0.
\end{equation}
Assume by contradiction that $\underline w(x_0)=0$, hence there exists a sequence $y_h\to x_0$ such that $\lim_{h\to 0}
w_{1,h}(y_h)=0$. By \eqref{Holder} with $u_h=w_{1,h}$ and $f=-\l_{1,h}w_{1,h}$ we get
\begin{align*}
|w_{1,h}(x_h)-w_{1,h}(y_h)|&\le  C \frac{|x_h-y_h|^\d}{R}\left(\max_{B_R} w_{1,h}+\frac{R}{\a_0}\left\{\sum_{x\in\Og_h}h^n|\l_{1,h}w_{1,h}|^{n}\right\}^{\frac{1}{n}}\right)\\
&\le C \frac{|x_h-y_h|^\d}{R}\left(1+\frac{R}{\a_0}|\l_{1,h}|\right).
\end{align*}
Since $\lim_{h\to 0}(w_{1,h}(x_h)-w_{1,h}(y_h))=1$ we get a contradiction for $h$ sufficiently small and
therefore \eqref{conv6}.\\
We are in a position to apply the maximum principle for the continuous problem (see \cite{BD}), and so
we obtain that $\underline w>0$.
But \eqref{conv3} and the positivity of $\underline w$ give a contradiction to  the definition of $\l_1$.
By \eqref{conv1} and  \eqref{conv2} we get $\lim_{h\to 0}\l_{1,h}=\l_1$. \\
By \eqref{Holder} and a local boundary estimate for $w_{1,h}$, see \cite[Thm. 5.1]{KT0} and \cite[Thm.3.]{KT1}, we get the equi-continuity of the family $\{w_{1,h}\}$ and therefore the uniform convergence, up to a subsequence,  of $w_{1,h}$ to $w_1$ with $\|w_{1}\|_\infty=1$.
The simplicity of the eigenfunction associated to the  principal eigenvalue $\l_1$ gives the uniform convergence of
 all the sequence $w_{1,h}$ to $w_1$.
\end{Proof}

\section{An algorithm for computing the principal eigenvalue}\label{sect4}
In this section we discuss an algorithm  for the computation of the principal eigenvalue based on  the inf-sup formula
\eqref{PEd_char2}. In fact  we show that this formula results in a finite dimensional nonlinear optimization problem.\par
\subsection{Discretization in one dimension.}
We first present the scheme in one dimension.
Note that  since the eigenfunction corresponding to the principal eigenvalue vanishes on the boundary of $\Og_h$ and it  is strictly positive inside,  then the  minimization in \eqref{PEd_char2}   can be restricted to  the internal points.
By the formula \eqref{finitediff} and  the homogeneity of $\cF$, we have
\[
   \frac{F_h[u](x_i)}{u(x_i)}={\cF}\left(x_i, 1, \frac{u(x_i+h)-u(x_i-h)}{2hu(x_i)},\frac{u(x_i+h)+u(x_i+h)}{h^2u(x_i)}-\frac{2}{h^2}\right).
\]
We identify the function $u(x)$ with the values $U_i$, $i=0,\dots, N_h+1$, at the points of the grid (with $U_0=U_{N_h+1}=0$).  Assume that $\cF(x,z,q,s)$ is  linear or more generally convex  in $(q,s)$. Then the functions $G:\R^{N_h}\to\R^{N_h}$, defined by
\[ G_i(x,U_1, \dots, U_{N_h})={\cF}\left(x_i, 1, \frac{U_{i+1}-U_{i-1}}{2hU_i},\frac{U_{i+1}+U_{i-1}}{h^2U_i}-\frac{2}{h^2}\right).\]
for $i=1,\dots,{N_h}$,  is  either  linear or respectively convex   in $U_{i+1}$, $U_{i-1}$. Moreover, since $U_i>0$, $G$ is also convex in $U_i$.
Taking the maximum of the functions $G_i$ over the internal nodes of the grid gives a convex function $\cG:\R^{N_h}\to\R$ defined by
\begin{equation}\label{minmax_2}
   \cG(U_1,\dots, U_{N_h}) =\max_{i=1,\dots, {N_h}}G_i(x_i, U_1, \dots, U_{N_h})
\end{equation}
Hence  the computation
of $\l_{1,h}$ is equivalent to the minimization of  the   convex function $\cG$ of ${N_h}$ variables:  this problem can be solved by means of standard algorithms in convex optimization. Note also that the minimum is unique and the map    is sparse, in the sense that
the value of $\cG$  at $U_i$ depends only on the values at $U_{i-1}$ and $U_{i+1}$.\\
In general, if $\cF(x,z,q,s)$ is   not convex, the computation of the principal eigenvalue is equivalent to
the solution of a min-max problem in $\R^{N_h}$.\\
 To solve min-max problem we use the routine \verb"fminmax" available  in the Optimization Toolbox of MATLAB. This routine is  implemented on a laptop and therefore the number of variables is modest. A better implementation of the minimization procedure which takes advantage of the sparse structure of the problem would allow to solve  larger problems.\vskip  4pt
\textbf{Example 1.}
To validate the algorithm we begin by studying the eigenvalue problem:
\[
    \left\{
    \begin{array}{ll}
    w^{\prime\prime}+\l w=0\quad& x\in (0,1), \\
    w(x)=0  & x=0,\,1.
    \end{array}
\right.
\]
In
this case the eigenvalue and the corresponding eigenfunction are  given by
\[\l_1= \pi^2,  \qquad w_1(x)= \sin(\pi x).\]
Note that since the eigenfunctions are  defined up to multiplicative constant, we normalize the value by taking $\|w_1\|_\infty=\|w_{1,h}\|_\infty=1$
(the constraint for $w_{1,h}$ is included in the routine \verb"fminmax").
Given a discretization step $h$ and the corresponding grid points $x_i=ih$, $i=0,\dots,N_h+1$, the minimization problem \eqref{minmax_2}  is
\[
   \l_{1,h}=-\min_{U\in\R^{N_h}}\left[\max_{i=1,\dots,N_h}\frac{U_{i+1}+U_{i-1}-2U_{i}}{h^2 U_i}\right]
\]
(with $U_0=U_{N_h+1}=0$).
In Table 4.1, we compare the exact solution with the approximate one  obtained by the scheme  \eqref{minmax_2}.
We report   the  approximation error  for  $\l_1$ and $w_1$ (in $L^\infty$-norm and $L^2$-norm) and the   order of convergence for
$\l_1$.  We can observe an order of convergence close to $2$ for $\l_1$ and therefore  equivalent to one obtained
 by  discretization  of  the  Rayleigh quotient via finite elements (see \cite{B}).
\begin{table}[ht!]
\begin{center}
\begin{tabular}{|c|c|c|c|c|}\hline
$h$ & $Err(\l_1)$ & $Order(\l_1)$  &$Err_\infty(w_1)$& $Err_2(w_1)$\\
\hline \hline
$1.00 \cdot 10^{-1}  $ &    $8.0908 \cdot 10^{-2}$   &             &   $3.3662\cdot 10^{-11}$ & $5.7732\cdot 10^{-11}$ \\ \hline
$5.00 \cdot 10^{-2}  $ &    $2.0277 \cdot 10^{-2}$ &1.9964 &  $1.4786\cdot 10^{-10} $&  $3.8119\cdot 10^{-10}$\\ \hline
$2.50  \cdot 10^{-2} $ &    $5.0723 \cdot 10^{-3}$ &1.9991 &  $6.6613\cdot 10^{-16}$&   $1.8731\cdot 10^{-15}$\\ \hline
$1.25  \cdot 10^{-2} $ &    $1.2683\cdot 10^{- 3}$ &1.9998 &  $1.5543\cdot 10^{-15}$ &  $6.2524\cdot 10^{-15}$\\ \hline
$6.25  \cdot 10^{-3} $ &    $3.1708\cdot 10^{-4}$  &1.9999 &  $1.2212 \cdot 10^{-15}$&  $7.1576\cdot 10^{-15}$\\ \hline
\end{tabular}
\label{tab1}
\end{center}
\caption{Space step (first column),  eigenvalue  error  (second column),   convergence  order (third column), eigenfunction error   in $L^\infty$  (fourth column), eigenfunction error   in $L^2$ (last column)}
\end{table}
\vskip  4pt
\textbf{Example 2. }
In this example we consider the eigenvalue problem for a linear equation with a discontinuous coefficient
\begin{equation}\label{eqdisc}
    \left\{
    \begin{array}{ll}
    a(x) w^{\prime\prime}+\l w=0\quad& x\in (0,\pi), \\
    w(x)=0  & x\in \{0,\,\pi\},
    \end{array}
\right.
\end{equation}
where
$$
a(x)=
\left\{\begin{array}{lc}
1 & \mbox{for }\ x\in [0,\frac{\pi}{2k}),\\
2 & \mbox{for }\ x\in [\frac{\pi}{2k},\pi],
\end{array}
\right.
$$
and $k:=\frac{2+\sqrt{2}}{2\sqrt{2}}>1$.
\begin{Proposition}
The principal eigenvalue $\l_1$ associated to problem \eqref{eqdisc} is  given by $k^2=\frac{3+2\sqrt{2}}{4}$.
\end{Proposition}
\begin{Proof}
Let $$w(x)=\left\{\begin{array}{lc}
\sin(kx) & \mbox{for }\ x\in [0,\frac{\pi}{2k}),\\
b\sin(\frac{kx}{\sqrt 2}+c) & \mbox{for }\ x\in [\frac{\pi}{2k},\pi].
\end{array}
\right.
$$
We choose  $b$ and $c$ such that  $w(0)=w(\pi)=0$ and $w$ continuous in
$\frac{\pi}{2k}$. Imposing these conditions we get
$$\frac{k\pi}{\sqrt 2}+c=\pi,\quad \mbox{and}\quad   b\sin(\frac{\pi}{2\sqrt 2}+c)=1,$$
i.e. $$c=\pi(1-\frac{k}{\sqrt 2})=\pi(\frac{2-\sqrt{2}}{4}).$$
Furthermore, using that $\frac{\pi}{2\sqrt{2}}+c=\frac{\pi}{2}$, we get
$$\lim_{x\rightarrow \frac{\pi}{2k}^-}w^\prime(x)=k\cos(\frac{\pi}{2})=0,\ \lim_{x\rightarrow \frac{\pi}{2k}^+}w^\prime(x)= bk\cos(\frac{\pi}{2\sqrt{2}}+c)=0,$$
hence  $w\in C^1([0,\pi])$.\\
On the other hand, since  $w$ is not $C^2$ in $\frac{\pi}{2k}$, we show that it satisfies \eqref{eqdisc} in the sense of viscosity solutions.
For any
$(p,q)\in J^{2,+}w(\frac{\pi}{2k})$, we get $p=0$ and $q\geq-\frac{k^2}{2}$.
This implies that for both $a=1$ and $a=2$:
$$aq\geq -k^2w(\frac{\pi}{2k}),$$ so $w$ is a subsolution.
For any  $(p,q)\in J^{2,-}w(\frac{\pi}{2k})$, we get $p=0$ and $q\leq-k^2$.
This implies that for both $a=1$ and $a=2$:
$$aq\leq -k^2w(\frac{\pi}{2k}),$$ so $w$ is a supersolution.
\end{Proof}
%
%
In Table 4.2, we compare the exact solution with the approximate one  obtained by means of the scheme
\[
   \l_{1,h}=-\min_{U\in\R^{N_h}}\left[\max_{i=1,\dots,N_h}a(ih)\frac{U_{i+1}+U_{i-1}-2U_{i}}{h^2 U_i}\right]
\]
(with $U_0=U_{N^h+1}=0$).  The rates are not very good, but  the problem
  is out of our setting since $F$ is discontinuous  and the error is very sensible to the chosen grid.
In Figure 4.1, we report the graph of the exact and   approximate eigenfunctions for $h=0.1$.
\begin{table}[ht!]
\begin{center}
\begin{tabular}{|c|c|c|c|c|}\hline
$h$ & $Err(\l_1)$ & $Order(\l_1)$     &$Err_\infty(w_1)$& $Err_2(w_1)$\\
\hline \hline
$ 0.1571  $   &$0.1197$   &    &  $0.0213$  &  $0.0563$ \\ \hline
$ 0.0785  $   &$ 0.0476$  & $1.3303$  &  $0.0090$  &  $0.0383$\\ \hline
$ 0.0393 $    &$0.0347$   & $0.4576$  &  $0.0065$  &  $0.0391$\\ \hline
$ 0.0196 $    &$ 0.0157$  & $1.1417$  &  $0.0030$  &  $0.0264$\\ \hline
$ 0.0098  $   &$0.0061$   & $1.3596$  &  $0.0012$  &  $0.0149$\\ \hline
\end{tabular}
\end{center}
\label{tab2}
\caption{Space step (first column), eigenvalue error   (second column),  eigenfunction error in $L^\infty$  (fourth column),  eigenfunction error  in $L^2$ (last column)}
\end{table}
\begin{figure}[ht!!]\label{fig1}
\begin{center}
\epsfig{figure=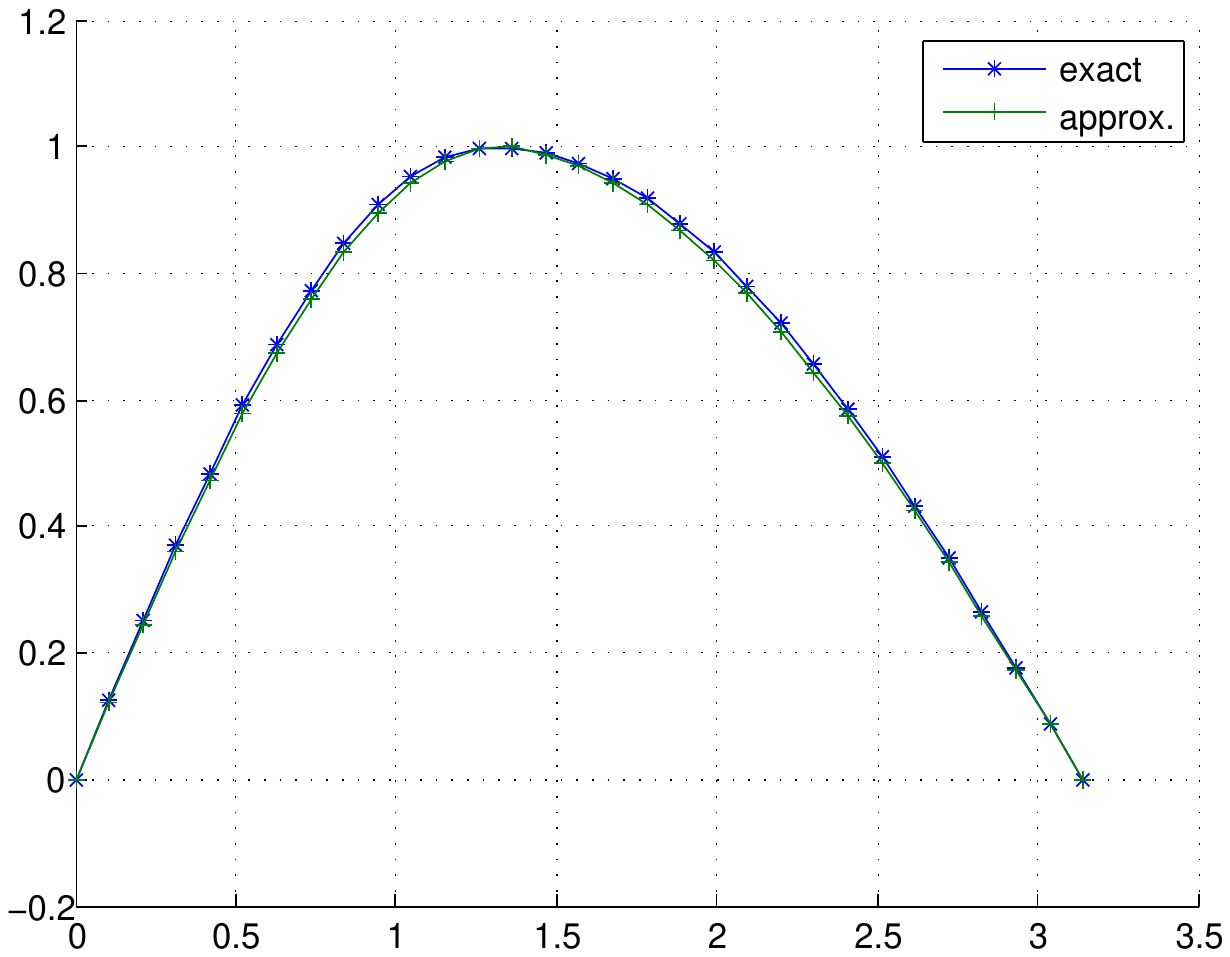,width=8cm}
\caption{ {{Exact and approximate eigenfunctions for  $h=10^{-1}$}}  }
\end{center}
\end{figure}
\vskip  4pt
\textbf {Example 3.}
The Fuc\^ik spectrum of $\Delta$  is the set of pairs $(\mu, \a\mu) \in\R^2$
for which the equation
\[
-\Delta u=\mu u^+-\a\mu u^-
\]
has a non-zero solution, where $u^+(x) = \max\{u(x), 0\}$ and $u^-(x) =
\max\{-u(x), 0\}$. For fixed$\a>0$ the previous problem is equivalent to \begin{align*}
    \min\{\Delta u,\frac 1 \a \Delta u\}+\mu u=0\quad &\text{if $\a\ge 1$},\\
    \max\{\Delta u,\frac 1 \a \Delta u\}+\mu u=0\quad &\text{if $\a\le1$}.
\end{align*}
For details see \cite{BEQ}.
 Hence the Fuc\^ik spectrum can be seen as the spectrum of a nonlinear operator  involving  the maximum or minimum of two linear operators. To find the corresponding  principal eigenvalue we apply the scheme \eqref{minmax_2}. In Table 4.3, we report the corresponding approximation error for $\l_1$ in the case $\a=1/2$ and $\Omega=[0,\pi]$ (by the convexity of the solution the eigenvalue for the continuous
problem coincides with the one of the second derivative operator in $[0,\pi]$ i.e. $\l_1=1$).
\begin{table}[ht!]
\begin{center}
\label{tab3}
\begin{tabular}{|c|c|c|}\hline
$h$ & $Err(\l_1)$ & $Order(\l_1)$  \\
\hline \hline
$1.00 \cdot 10^{-1}  $ &    $0.0809$  &                 \\ \hline
$5.00 \cdot 10^{-2}  $ &    $0.0203$ &  1.9964      \\ \hline
$2.50  \cdot 10^{-2} $ &    $0.0051$ &  1.9991       \\ \hline
$1.25  \cdot 10^{-2} $ &    $0.0013$ &  1.9998       \\ \hline
$6.25  \cdot 10^{-3} $ &    $0.0003$  & 2.0000    \\ \hline
\end{tabular}
\end{center}
\caption{Space step  (first column),  eigenvalue error  (second column),   convergence  order  (third column) for the Fuc\^ik spectrum with $\a=1/2$}
\end{table}

\vskip 4pt
\textbf{Example 4. }
Consider  the eigenvalue problem for the $p$-Laplace operator
\[
 \mathrm{div}(|Dw_1|^{p-2}Dw_1)+\lambda_1 |w_1|^{p-2}w_1=0.
\]
This example does not fit exactly in the framework of this paper since the operator is not uniformly elliptic.
However,
the following  formula
\[
     \l_{p,h}:=-\inf_{\phi>0}\sup_{y\in\Og_h}\left\{\frac{F_{h,p}[\phi](y)}{\phi(y)^{p-1}}\right\}
\]
where $F_{h,p}$ is a finite-difference approximations of $F_p$ produces a good approximation of the principal eigenvalue of the $p$-Laplace operator in the interval $(a,b)$ whose exact value is given by
 \[
 \sqrt[p]{\l_p} =\frac{2\pi\sqrt[p]{p-1}}{(b-a)p\sin(\frac{\pi}{p})}.
 \]
In Table 4.4 we report the approximation error and the corresponding order of convergence for the principal
 eigenvalue of the $p$-Laplace operator for $p=4$ (in this case $\l_4\approx 73.0568 $).\\
\begin{table}[ht!]
\begin{center}
\label{tab4}
\begin{tabular}{|c|c|c|}\hline
$h$ & $Err(\l_4)$ & $Order(\l_4)$  \\
\hline \hline
$1.00 \cdot 10^{-1}  $ &    $2.6770$  &                 \\ \hline
$5.00 \cdot 10^{-2}  $ &    $0.6210$ &  2.1079     \\ \hline
$2.50  \cdot 10^{-2} $ &    $0.1457$ &  2.0912       \\ \hline
$1.25  \cdot 10^{-2} $ &    $0.0347$ &  2.0724       \\ \hline
$6.25  \cdot 10^{-3} $ &    $0.0083$  & 2.0581    \\ \hline
\end{tabular}
\end{center}
\caption{Space step(first column), eigenvalue error   (second column),  convergence  order (third column) for the $p$-Laplace operator with $p=4$}
\end{table}
It is also known (see \cite{JPM}) that if $\Og$ is a ball, the eigenfunction $w_p$ corresponding to the eigenvalue $\l_p$ converges for $p\to \infty$ to $d(x,\partial\Og)$. In Figure 4.2, we draw   approximations of $w_p$ computed by the scheme  for various values of $p$ and we observe the convergence of these functions to $d(x,\{0,1\})$ for $p$ increasing, as expected by the theory.
\begin{figure}[ht!!]\label{fig2}
\begin{center}
\epsfig{figure=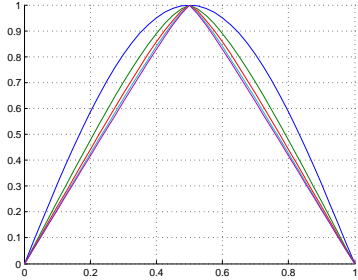,width=8cm}
\caption{ {{Approximate eigenfunction $u_{p,h}$ for $p=2,4,6,8,10$ and $h=10^{-3}$}}  }
\end{center}
\end{figure}
\subsection{Discretization in higher dimension.}
We now consider the eigenvalue problem in $\R^N$. Arguing as in the 1-dimensional case we write
  \begin{equation}\label{minmax_3}
 \small{\begin{split}
   \frac{F_h[u](x)}{u(x)} = {\cF}\Big(x, 1, &\left\{\frac{u(x+hy)-u(x-hy)}{2h|y|u(x)}\right\}_{y\in Y}, \left\{\frac{u(x+hy)+u(x-hy)}{h^2|y|^2u(x)}-\frac{2}{h^2}\right\}_{y\in Y}\Big)
   \end{split}}
  \end{equation}
for $i=1,\dots,N_h$  and    $F_h$ defined as  in \eqref{finitediff},  $Y$   the stencil and  ${N_h}$    the cardinality of $\Og_h$. Hence if  the function
$\cF(x,z,\{q_y\}_{y\in Y},\{s_y\}_{y\in Y})$ is   linear or more generally convex in the variables $q_y$ and $s_y$, $y\in Y$,
then the computation of the principal eigenvalue $\l_{1,h}$ is equivalent to the minimization with respect to  the vector $U\in\R^{N_h}$  of  the   convex function  $\cG:\R^{N_h}\to\R$ obtained by taking the maximum with respect to $x\in \Og_h$ in \eqref{minmax_3}.
Therefore   this problem can be solved by means of some standard algorithms in convex optimization.\vskip 4pt
\textbf{Example 5. }
Consider the problem
\[
    \left\{
    \begin{array}{ll}
   \Delta w+\l w=0\quad& x\in (0,1)^2, \\
    w(x)=0  & x\in \partial ((0,1)^2).
    \end{array}
\right.
\]
The eigenvalue and the corresponding eigenfunction   are  given by
\[\l_1=2 \pi^2,  \qquad w(x_1,x_2)= \sin(\pi x_1)\sin(\pi x_2)\]
(the eigenfunctions  are  normalized  by taking $\|w\|_\infty=\|w_{h}\|_\infty=1$). We use a standard five-point formula
for the discretization of the Laplacian.
In Table 4.5, we compare the exact solution with the approximate one  obtained by the scheme  \eqref{minmax_2}.
We report   the  approximation error  for  $\l_1$ and $w_1$ (in $L^\infty$-norm and $L^2$-norm) and the   order of convergence for
$\l_1$.  We can observe an order of convergence close to $2$ for $\l_1$ and therefore  equivalent to one obtained
 by  discretization  of  the  Rayleigh quotient via finite elements, see \cite{B}.
\begin{table}[ht!]
\begin{center}
\label{tab5}
\begin{tabular}{|c|c|c|c|c|}\hline
$h$ & $Err(\l_1)$ & $Order(\l_1)$  &$Err_\infty(w)$& $Err_2(w)$\\
\hline \hline
$2.00 \cdot 10^{-1}  $ &    $0.4469 $   &          &  $ 0.0801$ &  $0.2256$ \\ \hline
$1.00 \cdot 10^{-1}  $ &    $ 0.1338$   & $1.7397$ &  $ 0.0203$ &  $0.1137$\\ \hline
$5.00 \cdot 10^{-2}  $ &    $0.0368 $   & $1.8629$ &  $ 0.0056$ &  $0.0590$\\ \hline
$2.50  \cdot 10^{-2} $ &    $0.0097$    & $1.9297$ &  $ 0.0015$ &  $0.0301$\\ \hline
\end{tabular}
\end{center}
\caption{Space step  (first column),   eigenvalue   error (second column),   convergence  order (third column),   eigenfunction error in $L^\infty$  (fourth column),   eigenfunction error in $L^2$ (last column)}
\end{table}
\vskip 4pt
\textbf{Example 6. }
We consider the eigenvalue problem for the  Ornstein-Uhlenbeck operator
\[
\Delta w-x\cdot Dw+\l w=0, \qquad x\in (-1,1)^2
\]
 with homogeneous boundary conditions. The eigenvalue and the corresponding eigenfunction   are  given by
\[\l_1=4,  \qquad w(x_1,x_2)=(1-x_1^2)(1-x_2^2), \]
with the eigenfunctions  normalized  by taking $\|w\|_\infty=\|w_{1,h}\|_\infty=1$. The Laplacian
 is discretized by a five-point formula. In Table 4.6, we report the
 approximation error  for  $\l_1$ and the corresponding   order of convergence.
\begin{table}[ht!]
\begin{center}
\label{tab6}
\begin{tabular}{|c|c|c|c|c|}\hline
$h$ & $Err(\l_1)$ & $Order(\l_1)$  \\
\hline \hline
$4.00 \cdot 10^{-1}  $ &    $0.1524$   &                    \\ \hline
$2.00 \cdot 10^{-1}  $ &    $0.0392$   &  $1.9592$  \\ \hline
$1.00 \cdot 10^{-1}  $ &    $0.0103$   &  $1.9250$  \\ \hline
$5.00  \cdot 10^{-2} $ &    $0.0027$   &  $1.9580$  \\ \hline
\end{tabular}
\end{center}
\caption{Space step  (first column), eigenvalue  error  (second column),    convergence  order  (third column)}
\end{table}


\end{document}